\documentclass[twosided,reqno]{amsart}

\usepackage{amsmath}
\usepackage{amssymb}
\usepackage{amsthm}

\newcommand{\Cosh}{{\mathrm{cosh}}}
\newcommand{\Sinh}{{\mathrm{sinh}}}
\newcommand{\Sin}{{\mathrm{Sin}}}
\newcommand{\Cos}{{\mathrm{Cos}}}

\theoremstyle{theorem}
\newtheorem{theorem}{\scshape Theorem }[section]

\newtheorem{corollary}[theorem]{\scshape Corollary}

\theoremstyle{definition}

\numberwithin{equation}{section}

\begin{document}

\title[On the $q$-analogue of Laplace transform]{On the $q$-analogue of Laplace transform}

\author{Won Sang  Chung$^1$}
\address{$^1$ Department of Physics, Gyeongsang National University, Jinju 660,701, Republic of
Korea.}
\email{mimip4444@hanmail.net}

\author{Taekyun Kim$^2$}
\address{$^2$ Department of Mathematics, Kwangwoon University, Seoul 139-701, Republic of Korea.}
\email{tkkim@kw.ac.kr}


\maketitle

\begin{abstract}
In this paper, we consider a $q$-analogue of Laplace transform and we investigate some properties of $q$-Laplace transform. From our investigation, we derive some interesting formulae related to $q$-Laplace transform.
\end{abstract}

\section{Introduction}
For $q\in[0,1]$, we define the $q$-Shifted factorials by $(a:q)_0=1$, $(a:q)_n=\prod_{i=0} ^{n-1}(1-aq^i)$, and $(a:q)_{\infty}=\lim_{n\rightarrow \infty}(a:q)_n=\prod_{i=0} ^{\infty}(1-aq^i)$. If $x$ is classical object, such as a complex number, its $q$-version is defined by $[x]_q=\frac{1-q^x}{1-q}$. As is well-known, the {\it{$q$-exponential functions}} are given by
\begin{equation}\label{1}
E_q(-z)=\left((1-q)z:q)\right)_{\infty}=\sum_{n=0} ^{\infty}\frac{(-1)^nq^{\binom{n}{2}}}{[n]_q !}z^n,
\end{equation}
and
\begin{equation}\label{2}
e_q(z)=\frac{1}{\left((1-q)z:q\right)_{\infty}}=\sum_{n=0} ^{\infty}\frac{1}{[n]_q!}z^n,{\text{ (see [1-8])}},
\end{equation}
where $[n]_q!=[n]_q[n-1]_q\cdots[2]_q[1]_q$.

For $t,x,y\in{\mathbb{R}}$ and $n \in {\mathbb{Z}}\geq0$, the {\it{$q$-binomial theorem}} is also given by
\begin{equation}\label{3}
(x+y)_q ^n=\prod_{j=0} ^{n-1}\left(x+q^jy\right)=\sum_{k=0} ^n \binom{n}{k}_qx^{n-k}q^{\binom{2}{k}}y^k,
\end{equation}

\begin{equation}\label{4}
(1+x)_q ^t=\frac{(x:q)_{\infty}}{(q^tx:q)_{\infty}}, \text{ and } ~\frac{1}{(x-y)_q ^n}=\sum_{k=0} ^{\infty}\binom{n+k-1}{k}_qx^{n-k}y^k,
\end{equation}
where
\begin{equation*}
\binom{n}{k}_q=\frac{[n]_q !}{[n-k]_q ![k]_q !}=\frac{[n]_q[n-1]_q\cdots[n-k+1]_q}{[k]_q !},{\text{ (see [2,6])}}.
\end{equation*}
From \eqref{4}, we note that
\begin{equation}\label{5}
(1+x)_q ^{s+t}=(1+x)_q ^s(1+q^sx)_q ^t.
\end{equation}
For $s \in {\mathbb{C}}$ with $s>0$, the {\it{gamma function}} is defined by
\begin{equation}\label{6}
\Gamma(s)=\int_0 ^{\infty} e^{-t}t^{s-1}dt.
\end{equation}
By \eqref{6}, we get
\begin{equation*}
\Gamma(s+1)=s\Gamma(s),~\Gamma(n+1)=n!,~(n \in{\mathbb{N}}).
\end{equation*}
The {\it{Jackson $q$-derivative }} is defined by
\begin{equation}\label{7}
D_qf(x)=\frac{\partial_qf(x)}{\partial_qx}=\frac{f(x)-f(qx)}{(1-q)x},{\text{ (see [1-8])}}.
\end{equation}
Note that $\lim_{q\rightarrow 1}D_qf(x)=f^{'}(x)$. The {\it{definite Jackson $q$-integral}} is given by
\begin{equation}\label{8}
\int_0 ^x f(t)d_qt=(1-q)\sum_{a=0} ^{\infty}f\left(q^ax\right)xq^a,{\text{ (see [3-4])}}.
\end{equation}
Thus, by \eqref{7} and \eqref{8}, we get
\begin{equation}\label{9}
\int_0 ^x \frac{\partial_q}{\partial_q t}f(t)d_q t=f(x)-f(0).
\end{equation}
The improper $q$-integral of $f$ is given by
\begin{equation}\label{10}
\int_0 ^{\frac{\infty}{a}}f(x)d_qx=(1-q)\sum_{n\in{\mathbb{Z}}}\frac{q^n}{a}f\left(\frac{q^n}{a}\right).
\end{equation}
Thus, by \eqref{10}, we get
\begin{equation}\label{11}
\begin{split}
\int_x ^{\infty}f(y)d_qy&=\int_0 ^{x\cdot\infty}f(y)d_1y-\int_0 ^x f(y)d_qy \\
&=x(1-q)\sum_{n=0} ^{\infty}q^{-n}f\left(q^{-n}x\right).
\end{split}
\end{equation}
Let $f$ be a function defined for $t \geq0$. Then the integral
\begin{equation*}
{\mathcal{L}}(f(t))=\int_0 ^{\infty}e^{-st}f(t)dt,{\text{ (see [9])}}
\end{equation*}
is said to be the {\it{Laplace transform}} of $f$, provided the integral converges.

In this paper, we consider a $q$-analogue of Laplace transform, which is called by {\it{$q$-Laplace transform}}, and we investigate some properties of $q$-Laplace transform. From our investigation, we derive some interesting formulae related to $q$-Laplace transform.

\section{$q$-Laplace transforms}

Now we consider two types of the $q$-Laplace transform. The {\it{$q$-Laplace transform of the first kind}} is defined by
\begin{equation}\label{12}
F(s)=L_q(f(t))=\int_0 ^{\infty}E_q(-qst)f(t)d_q(t),~(s>0).
\end{equation}
Then, by \eqref{12}, we get
\begin{equation}\label{13}
L_q(\alpha f(t)+\beta g(t))=\alpha L_q(f(t))+\beta_qLg(t)),
\end{equation}
where $\alpha,\beta$ are constants.

The {\it{$q$-extension of gamma function}} is defined by
\begin{equation}\label{14}
\Gamma_q(t)=\int_0 ^{\infty}x^{t-1}E_q(-qx)d_qx,~(t>0),{\text{ (see [2,6])}}.
\end{equation}
Thus, by \eqref{14}, we get
\begin{equation}\label{15}
\Gamma_q(t+1)=[t]_q\Gamma_q(t),~\Gamma_q(n+1)=[n]_q !,{\text{ (see [2,6])}}.
\end{equation}
From \eqref{2}, we can derive
\begin{equation}\label{16}
e_q(iz)=\sum_{n=0} ^{\infty}\frac{i^nz^n}{[n]_q!}=\sum_{n=0} ^{\infty}\frac{(-1)^n}{[2n]_q!}z^{2n}+i\sum_{n=0} ^{\infty}\frac{(-1)^nz^{2n+1}}{[2n+1]_q!}.
\end{equation}
By \eqref{16}, we define {\it{$q$-cosine}} and {\it{$q$-sine function}} as follows:
\begin{equation}\label{17}
\cos_q(z)=\sum_{n=0} ^{\infty}\frac{(-1)^nz^{2n}}{[2n]_q!}=\frac{1}{2}\left(e_q(iz)+e_q(-iz)\right),
\end{equation}
and
\begin{equation}\label{18}
\sin_q(z)=\sum_{n=0} ^{\infty}\frac{(-1)^n}{[2n+1]_q!}=\frac{1}{2i}\left(e_q(iz)-e_q(-iz)\right),
\end{equation}
where $i=\sqrt{-1}$. From \eqref{12} and \eqref{15}, we note that
\begin{gather}
L_q(1)=\frac{1}{s}~(s>0),~L_q(t)=\frac{1}{s^2}~(s>0),~L_q\left(e_q(-3t)\right)=\frac{1}{s+3}~(s>-3),\notag \\
L_q\left(\sin_q2t\right)=\frac{2}{s^2+4},~L_q(1+5t)=\frac{1}{s}+\frac{5}{s^2}~(s>0),\cdots.\notag
\end{gather}
We state the generalization of some of the preceding examples by mean of the next theorem. From this point on we shall also refrain from stating any restrictions on $s$; it is understood that $s$ is sufficiently restricted to guarantee the convergence of the appropriate $q$-Laplace transform.

For $\alpha\in{\mathbb{R}}$ with $\alpha>-1$, we have
\begin{equation}\label{19}
\begin{split}
L_q\left(t^{\alpha}\right)&=\int_0 ^{\infty}E_q(-qst)t^{\alpha}d_qt=\frac{1}{s^{\alpha+1}}\int_0 ^{\infty}E_q(-qt)t^{\alpha}d_qt \\
&=\frac{1}{s^{\alpha+1}}\Gamma_q(\alpha+1).
\end{split}
\end{equation}
In particular, $\alpha=n\in{\mathbb{N}}$, by \eqref{19}, we get
\begin{equation}\label{20}
L_q\left(t^n\right)=\frac{1}{s^{n+1}}\Gamma_q(n+1)=\frac{[n]_q!}{s^{n+1}}.
\end{equation}
Let us take $\alpha=-\frac{1}{2}$ and $\alpha=\frac{1}{2}$. Then we see that
\begin{equation}\label{21}
L_q\left(t^{-\frac{1}{2}}\right)=\frac{1}{s^{\frac{1}{2}}}\Gamma_q\left(\frac{1}{2}\right)=\frac{1}{\sqrt{s}}\Gamma_q\left(\frac{1}{2}\right),
\end{equation}
and
\begin{equation}\label{22}
L_q\left(t^{\frac{1}{2}}\right)=\frac{1}{s^{\frac{3}{2}}}\Gamma_q\left(\frac{3}{2}\right)=\frac{\Gamma_q\left(\frac{1}{2}\right)}{[2]_{q^{\frac{1}{2}}}\sqrt{s^3}}.
\end{equation}
By \eqref{1} and \eqref{2}, we get
\begin{equation}\label{23}
\begin{split}
L_s\left(e_q(at)\right)&=\int_0 ^{\infty} E_q(-qst)e_q(at)d_qt \\
&=\sum_{n=0} ^{\infty}\frac{a^n}{[n]_q!}\int_0 ^{\infty}E_q(-qst)t^nd_qt \\
&=\sum_{n=0} ^{\infty}\frac{a^n}{[n]_q!}\frac{\Gamma_q(n+1)}{s^{n+1}}=\frac{1}{s}\sum_{n=0} ^{\infty}\left(\frac{a}{s}\right)^n \\
&=\frac{1}{s}\left(\frac{1}{1-\frac{a}{s}}\right)=\frac{1}{s-a},
\end{split}
\end{equation}
and
\begin{equation}\label{24}
\begin{split}
L_q\left(E_q(at)\right)=&\int_0 ^{\infty}E_q(-qst)E_q(at)d_qt \\
=&\sum_{n=0} ^{\infty}\frac{(-1)^nq^{\binom{n}{2}}a^n}{[n]_q!}\int_0 ^{\infty}E_q(-qst)t^nd_qt\\
=&\sum_{n=0} ^{\infty}\frac{(-1)^n}{[n]_q!}q^{\binom{n}{2}}a^n\frac{\Gamma_q(n+1)}{s^{n+1}}=\sum_{n=0} ^{\infty}(-1)^nq^{\binom{n}{2}}\frac{a^n}{s^{n+1}}.
\end{split}
\end{equation}
Therefore, by \eqref{19}-\eqref{24}, we obtain the following theorem.
\begin{theorem}\label{thm1}
For $\alpha\in{\mathbb{R}}$ with $\alpha>-1$, we have
\begin{equation*}
L_q\left(t^{\alpha}\right)=\frac{1}{s^{n+1}}\Gamma_q(\alpha+1).
\end{equation*}
In particular, $\alpha=n\in{\mathbb{N}}$, we get
\begin{equation*}
L_q\left(t^n\right)=\frac{[n]_q!}{s^{n+1}}.
\end{equation*}
Moreover,
\begin{equation*}
L_q\left(t^{-\frac{1}{2}}\right)=\frac{1}{\sqrt{s}}\Gamma_q\left(\frac{1}{2}\right),~L_q\left(t^{\frac{1}{2}}\right)=\frac{\Gamma_q\left(\frac{1}{2}\right)}{[2]_{q^{\frac{1}{2}}}\sqrt{s^3}},
\end{equation*}
and
\begin{equation*}
L_s\left(e_q(at)\right)=\frac{1}{s-a},~L_q(E_q(at))=\sum_{n=0} ^{\infty}(-1)^nq^{\binom{n}{2}}\frac{a^n}{s^{n+1}}.
\end{equation*}
\end{theorem}
From \eqref{17} and \eqref{18}, we have
\begin{equation}\label{25}
\begin{split}
L_q(\cos_qat)=&\sum_{n=0} ^{\infty}\frac{(-1)^na^{2n}}{[2]_q!}\int_0 ^{\infty}E_q(-qst)t^{2n}d_qt\\
=&\sum_{n=0} ^{\infty}\frac{(-1)^na^{2n}}{[2n]_q!}\frac{1}{s^{2n+1}}\Gamma_q(2n+1)\\
=&\frac{1}{s}\sum_{n=0} ^{\infty}(-1)^n\left(\frac{a}{s}\right)^{2n}=\frac{1}{s}\frac{1}{1+\left(\frac{a}{s}\right)^2}=\frac{s}{s^2+a^2},
\end{split}
\end{equation}
and
\begin{equation}\label{26}
\begin{split}
L_q(\sin_qat)=&\sum_{n=0} ^{\infty}\frac{(-1)^na^{2n+1}}{[2n+1]_q!}\int_0 ^{\infty}E_q(-qst)t^{2n+1}d_qt\\
=&\sum_{n=0} ^{\infty}\frac{(-1)^na^{2n+1}}{[2n+1]_q!}\frac{\Gamma_q(2n+2)}{s^{2n+2}}\\
=&\frac{1}{s}\sum_{n=0} ^{\infty}(-1)^n\left(\frac{a}{s}\right)^{2n+1}=\frac{1}{s}\frac{\frac{a}{s}}{1+\left(\frac{a}{s}\right)^2}=\frac{a}{s^2+a^2}.
\end{split}
\end{equation}
Let us define {\it{hyperbolic $q$-cosine}} and {\it{hyperbolic $q$-sine functions}} as follows:
\begin{equation}\label{27}
\Cosh_qt=\frac{e_q(t)+e_q(-t)}{e},~\Sinh_qt=\frac{e_q(t)-e_q(-t)}{2}.
\end{equation}
From \eqref{12}, \eqref{23} and \eqref{27}, we note that
\begin{equation}\label{28}
\begin{split}
L_q\left(\Cosh_qat\right)=&\frac{1}{2}\left\{L_q\left(e_q(at)\right)+L_q\left(e_q(-at)\right)\right\} \\
=&\frac{1}{2}\left\{\frac{1}{s-a}+\frac{1}{s+a}\right\}=\frac{s}{s^2-a^2},
\end{split}
\end{equation}
and
\begin{equation}\label{29}
\begin{split}
L_q\left(\Sinh_qat\right)=&\frac{1}{2}\left\{L_q\left(e_q(at)\right)-L_q\left(e_q(-at)\right)\right\} \\
=&\frac{1}{2}\left\{\frac{1}{s-a}-\frac{1}{s+a}\right\}=\frac{a}{s^2-a^2}.
\end{split}
\end{equation}
Therefore, by \eqref{25}-\eqref{29}, we obtain the following theorem.
\begin{theorem}\label{thm2}
(Transforms of $q$-trigonometric function)
\begin{gather}
L_q\left(\cos_qat\right)=\frac{s}{s^2+a^2},~L_q\left(\sin_qat\right)=\frac{a}{s^2+a^2} \notag \\
L_q\left(\Cosh_qat\right)=\frac{s}{s^2-a^2},~L_q\left(\Sinh_qat\right)=\frac{a}{s^2-a^2}.\notag
\end{gather}
\end{theorem}
Now, we observe that
\begin{equation}\label{30}
\begin{split}
E_q(-qst)=&\sum_{n=0} ^{\infty}\frac{(-1)^nq^{\binom{n}{2}}}{[n]_q!}q^ns^nt^n=\sum_{n=0} ^{\infty}\frac{(-1)^nq^{\binom{n}{2}}}{[n]_q!}\left([n]_q(q-1)+1\right)s^nt^n \\
=&(q-1)\sum_{n=1} ^{\infty}\frac{(-1)^nq^{\binom{n}{2}}}{[n-1]_q!}s^nt^n+\sum_{n=0} ^{\infty}\frac{(-1)^nq^{\binom{n}{2}}}{[n]_q!}s^nt^n \\
=&-(q-1)st\sum_{n=0} ^{\infty}\frac{(-1)^nq^{\binom{n}{2}}}{[n]_q!}q^ns^nt^n+\sum_{n=0} ^{\infty}\frac{(-1)^nq^{\binom{n}{2}}}{[n]_q!}s^nt^n \\
=&-(q-1)stE_q(-qst)+E_q(-st).
\end{split}
\end{equation}
Thus, by \eqref{30}, we get
\begin{equation}\label{31}
\begin{split}
E_q(-qst)=&\frac{1}{1+(q-1)st}E_q(-st)=\frac{1}{(1+(q-1)st)(1+(q-1)q^{-1}st)}E_q\left(-q^{-1}st\right) \\
=&\cdots=\frac{1}{(1+(q-1)st)(1+(q-1)q^{-1}st)\cdots (1+(q-1)q^{-k}st)}E_q\left(-q^{-k}st\right).
\end{split}
\end{equation}
Now, we dicuss the $q$-differential equation. The main purpose of $q$-Laplace transform is in converting $q$-differential equation into simpler form which may be solved more easily. Like the ordinary Laplace transform, we can compute the $q$-Laplace transform of derivative by using the definition of $q$-Laplace transform.

Now, we oberve that
\begin{equation}\label{32}
\begin{split}
D_q\left(f(t)g(t)\right)=&\frac{\partial}{\partial_qt}(f(t)g(t))=\frac{\partial f(t)}{\partial_qt}g(t)+f(qt)\frac{\partial g(t)}{\partial_qt} \\
=&\left(D_qf(t)\right)g(t)+f(qt)\left(D_qg(t)\right).
\end{split}
\end{equation}
Thus, by \eqref{32}, we get
\begin{equation}\label{33}
\int_0 ^{\infty}\left(D_qf(t)\right)g(t)d_qt=\int_0 ^{\infty}D_q(f(t)g(t))d_qt-\int_0 ^{\infty} f(qt)\left(D_qg(t)\right)d_qt.
\end{equation}
It is easy to show that
\begin{equation}\label{34}
\begin{split}
D_qE_q(-qst)=&\frac{\partial}{\partial_qt}\sum_{n=1} ^{\infty}\frac{q^{\binom{n}{2}}(-1)^n}{[n]_q!}q^ns^nt^n=\sum_{n=1} ^{\infty}\frac{q^{\binom{n}{2}}(-1)^n}{[n-1]_q!}q^ns^nt^{n-1}\\
=&-s\sum_{n=0} ^{\infty}\frac{q^{\binom{n}{2}}(-1)^n}{[n]_q!}q^{n+1}s^nt^n.
\end{split}
\end{equation}
From \eqref{31}, \eqref{33} and \eqref{34}, we have
\begin{equation}\label{35}
\begin{split}
L_q\left(D_qf(t)\right)=&\int_0 ^{\infty}\left(D_qf(t)\right)E_q(-qst)d_qt=-f(0)-\int_0 ^{\infty}f(qt)\frac{\partial}{\partial_qt}E_q(-qst)d_qt \\
=&-f(0)+s\sum_{n=0} ^{\infty}\frac{q^{\binom{n+1}{2}}}{[n]_q!}(-1)^n\int_0 ^{\infty}f(qt)q^{n+1}s^nt^nd_qt\\
=&-f(0)+s\int_0 ^{\infty}f(t)\sum_{n=0} ^{\infty}\frac{(-1)^nq^{\binom{n}{2}}}{[n]_q!}(qst)^nd_qt \\
=&-f(0)+s\int_0 ^{\infty}f(t)E_q(-qst)d_qt=-f(0)+sL_q(f(t)).
\end{split}
\end{equation}
If we replace $f(t)$ by $D_qf(t)$, we see that
\begin{equation}\label{35}
\begin{split}
L_q\left(D_q ^2f(t)\right)=&-\left(D_qf\right)(0)+s\int_0 ^{\infty}\left(D_qf(t)\right)E_q(-q(t))d_qt\\
=&-\left(D_qf\right)(0)+sL_q\left(D_qf(t)\right)\\
=&-\left(D_qf\right)(0)+s(-f(0))+sL_q(f(t)) \\
=&-\left(D_qf\right)(0)-sf(0)+s^2L_q(f(t)).
\end{split}
\end{equation}
Continuing this process, we get
\begin{equation}\label{37}
L_q\left(f^{(n)}(t)\right)=s^nL_q\left(f(t)\right)-\sum_{i=0} ^{n-1}s^{n-1-i}f^{(i)}(0),
\end{equation}
where $f^{(n)}(t)=\left(\frac{\partial_q}{\partial_qt}\right)^nf(t)=D_q ^nf(t)$, $f^{(n)}(0)=\left.f^{(n)}(t)\right| _{t=0}$.

A function $f$ is said to be {\it{of exponential order $c$}} if there exists $c$, $M>0$ and $T>0$ such that
\begin{equation*}
|f(t)| \leq Me^{ct}{\text{ for all }}t>T.
\end{equation*}
If $f(t)$ is piecewise continuous on the interval $(0,\infty)$ and of exponential order $c$, then $L_q(f(t))$ exists for $s>c$. Therefore, by \eqref{37}, we obtain the following theorem.
\begin{theorem}\label{thm3}
If $f,~f^{'},\cdots,f^{(n-1)}$ are continuous on $(0,\infty)$ and are of exponential order and if $f^{(n)}(t)$ is piecewise continunus on $(0,\infty)$, then we have
\begin{equation*}
L_q\left(f^{(n)}(t)\right)=s^nL_q\left(f(t)\right)-\sum_{i=0} ^{n-1}s^{n-1-i}f^{(i)}(0),
\end{equation*}
where $f^{(n)}(t)=\left(\frac{\partial_q}{\partial_qt}\right)^nf(t)$.
\end{theorem}
Let us consider the following $q$-derivative in $s$:
\begin{equation}\label{38}
\begin{split}
\frac{\partial F(q^{-1}s)}{\partial_q s}=&\int_0 ^{\infty}\left(\frac{\partial_q}{\partial_qs}E_q(-st)\right)f(t)d_qt \\
=&-\int_0 ^{\infty}tE_q(-qst)f(t)d_qt\\
=&-L_q\left(tf(t)\right),
\end{split}
\end{equation}
and
\begin{equation}\label{39}
\begin{split}
q\left(\frac{\partial _q}{\partial_q s}\right)^2F(q^{-2}s)=&q\left(\frac{\partial_q}{\partial_q s}\right)^2\int_0 ^{\infty}E_q\left(-q^{-1}st\right)f(t)d_qt\\
=&q\int_0 ^{\infty}\left(\left(\frac{\partial_q}{\partial_qs}\right)^2E_q\left(-q^{-1}st\right)\right)f(t)d_qt \\
=&q\sum_{n=0} ^{\infty}\frac{(-1)^{n+2}}{[n]_q!}q^{\binom{n}{2}}q^{n-1}s^n\int_0 ^{\infty}t^{n+2}f(t)d_qt\\
=&\int_0 ^{\infty} (-1)^2\left(\sum_{n=0} ^{\infty}\frac{(-1)^nq^{\binom{n}{2}}}{[n]_q!}q^ns^nt^n\right)t^2f(t)d_qt\\
=&(-1)^2\int_0 ^{\infty}E_q(-qst)t^2f(t)d_qt=(-1)^2L_q\left(t^2f(t)\right).
\end{split}
\end{equation}
Continuing this process, we get
\begin{equation}\label{40}
L_q\left(t^nf(t)\right)=(-1)^nq^{\binom{n}{2}}\left(\frac{\partial_q}{\partial_qs}\right)^nF\left(q^{-ns}\right).
\end{equation}
Therefore, by \eqref{40}, we obtain the following theorem.
\begin{theorem}[$q$-Derivative of Transforms]\label{thm4}
For $n \in {\mathbb{N}}$, we have
\begin{equation*}
L_q\left(t^nf(t)\right)=(-1)^nq^{\binom{n}{2}}\left(\frac{\partial_q}{\partial_qs}\right)^nF\left(q^{-n}s\right).
\end{equation*}
\end{theorem}
From Theorem \ref{thm4}, we note that
\begin{equation}\label{41}
\begin{split}
L_q\left(t^ne_q(at)\right)=&(-1)^nq^{\binom{n}{2}}\left(\frac{\partial_q}{\partial_qs}\right)^n
\int^{\infty}_{0}E_q(-q^{-n+1}st)e_q(at)d_qt  \\
=&(-1)^nq^{\binom{n}{2}}\left(\frac{\partial_q}{\partial_qs}\right)^n\left(\frac{1}{q^{-n}s-a}\right)\\
=&\frac{(-1)^nq^{\binom{n}{2}}[n]_q!(-1)^nq^{-n^2}}{(s-a)(q^{-1}s-a)\cdots(q^{-n}s-a)} \\
=&\frac{q^{-\binom{n+1}{2}}[n]_q!}{(s-a)(q^{-1}s-a)\cdots(q^{-n}s-a)},
\end{split}
\end{equation}
and
\begin{equation}\label{42}
\begin{split}
L_q\left(e_q(at)f(t)\right)=&\sum_{n=0} ^{\infty}\frac{a^n}{[n]_q!}L_q\left(t^nf(t)\right)\\
=&\sum_{n=0} ^{\infty}\frac{a^n}{[n]_q!}(-1)^nq^{\binom{n}{2}}\left(\frac{\partial_q}{\partial_qs}\right)^nF(q^{-n}s)\\
=&\sum_{n=0} ^{\infty}\frac{(-a)^n}{[n]_q!}q^{\binom{n}{2}}\left(\frac{\partial_q}{\partial_qs}\right)^nF(q^{-n}s).
\end{split}
\end{equation}
Therefore, by \eqref{41} and \eqref{42}, we obtain the following corollary.
\begin{corollary}\label{coro5}
For $n \in {\mathbb{N}}$, we have
\begin{equation*}
L_q\left(e_q(at)f(t)\right)=\sum_{n=0} ^{\infty}\frac{(-a)^n}{[n]_q!}q^{\binom{n}{2}}\left(\frac{\partial_q}{\partial_qs}\right)^nF(q^{-n}s).
\end{equation*}
In particular, for $f(t)=t^n~(n\in{\mathbb{N}})$, we have
\begin{equation*}
L_q\left(t^ne_q(at)\right)=\frac{q^{-\binom{n+1}{2}}[n]_q!}{(s-a)(q^{-1}s-a)\cdots(q^{-n}s-a)}.
\end{equation*}
\end{corollary}
For $t,s>0$, the {\it{$q$-beta function}} is defined by
\begin{equation}\label{43}
B_q(t,s)=\int_0 ^1 x^{t-1}(1-qx)_q ^{s-1} d_q x,{\text{ (see [2,6])}}.
\end{equation}
From \eqref{43}, we note that
\begin{equation}\label{44}
B_q(t,s)=\frac{\Gamma_q(t)\Gamma_q(s)}{\Gamma_q(t+s)},{\text{ (see [2,6])}}.
\end{equation}
It is easy to show that
\begin{equation}\label{45}
T_q ^{\alpha}(f)=\int_0 ^t (t-qs)_q ^{\alpha-1}f(s)d_qs=t\int_0 ^1 (t-qrt)_q ^{\alpha-1}f(rt)d_qr,~(\alpha>0).
\end{equation}
In particular, if we take $f(t)=t^{\beta}~(\beta>0)$, then
\begin{equation}\label{46}
\begin{split}
\int_0 ^t (t-qs)_q ^{\alpha-1}t^{\beta}d_qs=&t^{\alpha+\beta}\int_0 ^1 \left(1-qr\right)_q ^{\alpha-1}r^{\beta}d_qr \\
=&B_q(\alpha,\beta+1)t^{\alpha+\beta}.
\end{split}
\end{equation}
If functions $f$ and $g$ are piecewise continuous on $(0,\infty)$, then a {\it{special product}}, denote by $f{*}g$ is defined by the integral
\begin{equation}\label{47}
(f{*}g)(t)=\int_0 ^t f(\tau)g(t-\tau)d\tau,
\end{equation}
and is called {\it{convolution}} of $f$ and $g$.

Now, we consider the $q$-analogue of convolution of $f$ and $g$.

Let $f_1(t)=t^{\alpha}$, $g(t)=t^{\beta-1}$ $(\alpha, \beta>0)$. Then we define the {\it{$q$-convolution}} of $f$ and $g$ as follows:
\begin{equation}\label{48}
(f_1{*}g)(t)=\int_0 ^tf_1(\tau)g(t-q\tau)d_q\tau,
\end{equation}
where $g(t-q\tau)=(t-q\tau)_q ^{\beta-1}$. From \eqref{48}, we have
\begin{equation}\label{49}
\begin{split}
(f_1{*}g)(t)=&\int_0 ^t \tau^{\alpha}(t-\tau q)_q ^{\beta-1}d_q \tau=t\int_0 ^1 (rt)^{\alpha}(t-qrt)_q ^{\beta-1}d_q r \\
=&t^{\alpha+\beta}\int_0 ^1 r^{\alpha}(1-qr)^{\beta-1}d_qr=t^{\alpha+\beta}B_q(\alpha+1,\beta) \\
=&\frac{\Gamma_q(\alpha+1)\Gamma_q(\beta)}{\Gamma_q(\alpha+\beta+1)}t^{\alpha+\beta}.
\end{split}
\end{equation}
Thus, by \eqref{49}, we get
\begin{equation}\label{50}
\begin{split}
L_q(f_1{*}g)=&B_q(\alpha+1,\beta)\int_0 ^{\infty} E_q(-qst)t^{\alpha+\beta}d_qt\\
=&B_q(\alpha+1,\beta)\frac{1}{s^{\alpha+\beta+1}}\Gamma_q(\alpha+\beta+1)\\
=&\left(\frac{\Gamma_q(\alpha+1)}{s^{\alpha+1}}\right)\left(\frac{\Gamma_q(\beta)}{s^{\beta}}\right).
\end{split}
\end{equation}
By \eqref{19}, we see that
\begin{equation}\label{51}
L_q(f_1)=\frac{\Gamma_q(\alpha+1)}{s^{\alpha+1}},~L_q(g)=\frac{1}{s^{\beta}}\Gamma_q(\beta).
\end{equation}
Hence, by \eqref{50} and \eqref{51}, we get
\begin{equation}\label{52}
L_q(f_1{*}g)=L_q(f_1)L_q(g).
\end{equation}
Assume that $f(t)$ is of the type such that equation \eqref{52}. Then we have
\begin{equation}\label{53}
\begin{split}
f{*}g=&\int_0 ^t (t-\tau q)_q ^{\beta-1}f(\tau)d_q\tau=t\int_0 ^1 (t-qrt)_q^{\beta-1}f(rt)d_qr \\
=&T_q^{\beta}(f(t)),
\end{split}
\end{equation}
where $g(t)=t^{\beta-1}$, and
\begin{equation}\label{54}
\begin{split}
L_q\left(T_q ^{\beta} f(t)\right)=&L_q(f{*}g)=L_q(f)L_q(g)\\
=&\frac{\Gamma_q(\beta)}{s^{\beta}}L_q(f).
\end{split}
\end{equation}
If $f(t)=\sum_i a_i t^{\alpha_i}$, then we have
\begin{equation}\label{55}
\begin{split}
L_q(f{*}g)=&\sum_i a_iL_q\left(t^{\alpha_i}{*}g\right) \\
=&\sum_i a_i L_q\left(t^{\alpha_i}\right)L_q(g)\\
=&L_q\left(\sum_i a_i t^{\alpha_i}\right)L_q(g)\\
=&L_q(f)L_q(g).
\end{split}
\end{equation}
Therefore, by \eqref{55}, we obtain the following theorem.
\begin{theorem}\label{thm6}
For $f(t)=\sum_i a_i t^{\alpha_i}$ and $g(t)=t^{\beta-1}$, we have
\begin{equation*}
L_q(f{*}g)=L_q(f)L_q(g).
\end{equation*}
\end{theorem}
For $\beta=1$, let $T_q f(t)=T_q ^1f(t)$. Then we see that
\begin{equation}\label{56}
\begin{split}
L_q\left(T_q\sin_qt\right)=&L_q(\sin_qt{*}1)=L_q(\sin_qt)L_q(1)\\
=&\left(\frac{1}{s^2+1}\right)\left(\frac{1}{s}\right)=\frac{1}{s(s^2+1)},
\end{split}
\end{equation}
and
\begin{equation}\label{57}
L_q\left(T_q\sin_qt\right)=L_q\left(\int_0 ^t \sin_q\tau d_q\tau\right).
\end{equation}
By \eqref{56} and \eqref{57}, we get
\begin{equation}\label{58}
L_q\left(\int_0 ^t \sin_q\tau d_q\tau\right)=\frac{1}{s(s^2+1)}.
\end{equation}
Note that
\begin{equation}\label{59}
\begin{split}
&L_q\left(\int_0 ^t \left(1-\cos_q\tau\right)d_q \tau\right)=L_q\left(T_q\left(1-\cos_qt\right)\right) \\
=&L_q\left(\left(1-\cos_qt\right){*}1\right)=L_q\left(1-\cos_qt\right)L_q(1) \\
=&\left(\frac{1}{s}-\frac{s}{s^2+1}\right)\frac{1}{s}=\frac{1}{s^2(s^2+1)},
\end{split}
\end{equation}
and
\begin{equation}\label{60}
\begin{split}
&L_q\left(\int_0 ^t \left(\tau-\sin_qt\right)d_q\tau\right) \\
=&L_q\left(T_q\left(t-\sin_qt\right)\right)=L_q\left(\left(t-\sin_qt\right){*}1\right)=L_q\left(t-\sin_qt\right)L_q(1)\\
=&\left(\frac{1}{s^2}-\frac{1}{s^2+1}\right)\frac{1}{s}=\frac{1}{s^2(s^2+1)}\frac{1}{s}=\frac{1}{s^3(s^2+1)}.
\end{split}
\end{equation}
Therefore, by \eqref{59} and \eqref{60}, we obtain the following corollary.
\begin{corollary}\label{coro7}
For $t\geq 0$, we have
\begin{equation*}
L_q\left(\int_0 ^t \left(1-\cos_q\tau\right) d_q\tau\right)=\frac{1}{s^2(s^2+1)},
\end{equation*}
and
\begin{equation*}
L_q\left(\int_0 ^t \left(\tau-\sin_q\tau\right) d_q\tau\right)=\frac{1}{s^3(s^2+1)},
\end{equation*}
\end{corollary}
Let $u(t-a)$ be {\it{Heaviside function}} which defined by
\begin{equation}\label{61}
u(t-a)=\left\{
\begin{array}{cl}
1 &{\text{ if }}t\geq a, \\
0 &{\text{ if }}0 \leq t <a.
\end{array}\right.
\end{equation}
Then, we have
\begin{equation}\label{62}
\begin{split}
L_q(u(t-a))=&\int_0 ^{\infty}E_q(-qst)u(t-a)d_qt=\int_a ^{\infty}E_q(-qst)d_qt \\
=&\int_0 ^{\infty}E_q(-qst)d_qt-\int_0 ^a E_q(-qst)d_qt \\
=&\frac{1}{s}-\int_0 ^a E_q(-qst)d_qt=\frac{1}{s}-\sum_{n=0} ^{\infty}\frac{(-1)^nq^{\binom{n}{2}}}{[n]_q!}(qs)^n\int_0 ^a t^nd_qt\\
=&\frac{1}{s}+\frac{1}{s}\sum_{n=1} ^{\infty}\frac{(-1)^n}{[n]_q!}q^{\binom{n}{2}}s^na^n=\frac{1}{s}\sum_{n=0} ^{\infty}\frac{(-1)^nq^{\binom{n}{2}}}{[n]_q!}(as)^n\\
=&\frac{1}{s}E_q(-qas).
\end{split}
\end{equation}
Therefore, by \eqref{62}, we obtain the following theorem.
\begin{theorem}\label{thm8}
Let $u(t-a)$ be Heaviside function. Then we have
\begin{equation*}
L_q\left(u(t-a)\right)=\frac{1}{s}E_q(-qas).
\end{equation*}
\end{theorem}

\section{$q$-Laplace transform of the second kind}

In this section, we consider the $q$-Laplace transform with $e_q(-st)$ and is called $q$-Laplace transform of the second kind. The {\it{$q$-Laplace transform of the second kind}} is defined by
\begin{equation}\label{63}
{\tilde{F}}(s)={\tilde{L}}_q\left(f(t)\right)=\int_0 ^{\infty} e_q(-st)f(t)d_qt,~(s>0).
\end{equation}
As is well known, the {\it{$q$-gamma function of the second kind}} is defined by
\begin{equation}\label{64}
\gamma_q(t)=\int_0 ^{\infty}x^{t-1}e_q(-x)d_qx,~(t>0),{\text{ (see [2,6])}}.
\end{equation}
Thus, by \eqref{64}, we get
\begin{equation}\label{65}
\gamma_q(1)=1,~\gamma_q(t+1)=q^{-t}[t]_q\gamma_q(t),~\gamma_q(n)=q^{-\binom{n}{2}}\Gamma_q(n),~(n\in{\mathbb{N}}).
\end{equation}
For $\alpha\in{\mathbb{R}}$ with $\alpha>-1$, we have
\begin{equation}\label{66}
\begin{split}
{\tilde{L}}_q\left(t^{\alpha}\right)=&\int_0 ^{\infty}e_q(-st)t^{\alpha}d_qt=\frac{1}{s^{\alpha+1}}\int_0 ^{\infty} e_q(-t)t^{\alpha}d_qt \\
=&\frac{1}{s^{\alpha+1}}\gamma_q(\alpha+1).
\end{split}
\end{equation}
In particular, for $\alpha=n\in{\mathbb{N}}$, by \eqref{65} and \eqref{66}, we get
\begin{equation}\label{67}
\begin{split}
{\tilde{L}}_q\left(t^n\right)=&\frac{1}{s^{n+1}}\gamma_q(n+1)=\frac{1}{s^{n+1}}q^{-\binom{n+1}{2}}\Gamma_q(n+1)\\
=&\frac{1}{s^{n+1}}q^{-\frac{(n+1)n}{2}}[n]_q!.
\end{split}
\end{equation}
Let us take $\alpha=-\frac{1}{2}$. Then we have
\begin{equation}\label{68}
{\tilde{L}}_q\left(t^{-\frac{1}{2}}\right)=\frac{1}{\sqrt{s}}\gamma_q\left(\frac{1}{2}\right),~{\tilde{L}}_q\left(t^{\frac{1}{2}}\right)=\frac{1}{{\sqrt{s^3}}q[2]_{q^{\frac{1}{2}}}}\gamma_q\left(\frac{1}{2}\right).
\end{equation}
In a similar way, we can give the $q$-Laplace transform of the second kind for $e_q(at)$ and $E_q(at)$.
\begin{equation}\label{70}
\begin{split}
{\tilde{L}}_q\left(e_q(at)\right)=&\sum_{n=0} ^{\infty}\frac{a^n}{[n]_q!}\int_0 ^{\infty}e_q(-st)t^nd_qt=\sum_{n=0} ^{\infty}\frac{a^n}{[n]_q!}{\tilde{L}}_q(t^n)\\
=&\sum_{n=0} ^{\infty}\frac{a^n}{[n]_q!}\frac{q^{-\binom{n+1}{2}}}{s^{n+1}}\Gamma_q(n+1)=\sum_{n=0} ^{\infty}\frac{a^n}{s^{n+1}}q^{-\binom{n+1}{2}},
\end{split}
\end{equation}
and
\begin{equation}\label{71}
\begin{split}
{\tilde{L}}_q\left(E_q(at)\right)=&\sum_{n=0} ^{\infty}\frac{q^{\binom{n}{2}}a^n}{[n]_q!}\int_0 ^{\infty}t^ne_q(-st)d_qt=\sum_{n=0} ^{\infty}\frac{q^{\binom{n}{2}}a^n}{[n]_q!}{\tilde{L}}_q(t^n)\\
=&\sum_{n=0} ^{\infty}\frac{q^{\binom{n}{2}}a^n}{[n]_q!}\frac{q^{-\binom{n+1}{2}}}{s^{n+1}}\Gamma_q(n+1)=\frac{1}{s}\sum_{n=0} ^{\infty}\left(\frac{a}{qs}\right)^n\\
=&\frac{1}{s}\left(\frac{qs}{qs-a}\right)=\frac{q}{qs-a}.
\end{split}
\end{equation}
Therefore, by \eqref{66}-\eqref{71}, we obtain the following theorem.
\begin{theorem}
For $\alpha\in{\mathbb{R}}$ with $\alpha>-1$, we have
\begin{equation*}
{\tilde{L}}_q\left(t^{\alpha}\right)=\frac{1}{s^{\alpha+1}}\gamma_q(\alpha+1),~(s>0).
\end{equation*}
In particular, if $\alpha=n\in{\mathbb{N}}$, then we see that
\begin{equation*}
{\tilde{L}}_q\left(t^n\right)=\frac{1}{s^{n+1}}q^{-\binom{n+1}{2}}[n]_q!.
\end{equation*}
Moreover,
\begin{equation*}
{\tilde{L}}_q\left(e_q(at)\right)=\sum_{n=0} ^{\infty}\frac{a^n}{s^{n+1}}q^{-\binom{n+1}{2}},~{\tilde{L}}_q\left(E_q(at)\right)=\frac{q}{qs-a},
\end{equation*}
and
\begin{equation*}
{\tilde{L}}_q\left(t^{-\frac{1}{2}}\right)=\frac{1}{\sqrt{s}}\gamma_q\left(\frac{1}{2}\right),~{\tilde{L}}_q\left(t^{\frac{1}{2}}\right)=\frac{1}{{\sqrt{s^3}}q[2]_{q^{\frac{1}{2}}}}\gamma_q\left(\frac{1}{2}\right).
\end{equation*}
\end{theorem}
From \eqref{1}, we note that
\begin{equation}\label{72}
E_q(it)=\sum_{n=0} ^{\infty}\frac{q^{\binom{n}{2}}}{[n]_q!}i^nt^n.
\end{equation}
Thus, by \eqref{72}, we get
\begin{equation}\label{73}
E_q(it)=\sum_{n=0} ^{\infty}\frac{(-1)^nq^{\binom{2n}{2}}}{[2n]_q!}t^{2n}+i\sum_{n=0} ^{\infty}\frac{(-1)^nq^{\binom{2n+1}{2}}}{[2n+1]_q!}t^{2n+1}.
\end{equation}
By \eqref{73}, we can define new {\it{$q$-Sine}} and {\it{$q$-Cosine}} function as follows:
\begin{equation}\label{74}
\Cos_q(at)=\sum_{n=0} ^{\infty}\frac{q^{\binom{2n}{2}}(-1)^na^{2n}}{[2n]_q!}t^{2n}=\frac{1}{2}\left[E_q(iat)+E_q(-iat)\right],
\end{equation}
and
\begin{equation}\label{75}
\Sin_q(at)=\sum_{n=0} ^{\infty}\frac{q^{\binom{2n+1}{2}}(-1)^na^{2n+1}}{[2n+1]_q!}t^{2n+1}=\frac{1}{2}\left[E_q(iat)-E_q(-iat)\right].
\end{equation}
Now, we give the $q$-Laplace transform of the second kind for $q$-Sine and $q$-Cosine functions:

From \eqref{27} and \eqref{70}, we note that
\begin{equation}\label{76}
\begin{split}
{\tilde{L}}_q\left(\cos_qat\right)=&\frac{1}{2}\left[{\tilde{L}}_q\left(e_q(iat)\right)+{\tilde{L}}_q(-iat)\right]\\
=&\frac{1}{2}\left[\sum_{I=0} ^{\infty}\frac{q^{-\binom{n+1}{2}}}{s^{n+1}}(ia)^n+\sum_{n=0} ^{\infty}\frac{q^{-\binom{n+1}{2}}}{s^{n+1}}(-ia)^n\right] \\
=&\frac{1}{2}\times 2\sum_{n=0} ^{\infty}\frac{q^{-n(2n+1)}}{s^{2n+1}}\left(\frac{a}{s}\right)^{2n}=\frac{1}{s}\sum_{n=0} ^{\infty}(-1)^nq^{-n(2n+1)}\left(\frac{a}{s}\right)^{2n},
\end{split}
\end{equation}
and
\begin{equation}\label{77}
\begin{split}
{\tilde{L}}_q\left(\sin_qat\right)=&\sum_{n=0} ^{\infty}\frac{(-1)^na^{2n+1}}{[2n+1]_q!}\int_0 ^{\infty} e_q(-st)t^{2n+1}d_qt=\sum_{n=0} ^{\infty}\frac{(-1)^na^{2n+1}}{[2n+1]_q!}{\tilde{L}}_q(t^{2n+1})\\
=&\sum_{n=0} ^{\infty} \frac{(-1)^na^{2n+1}}{[2n+1]_q!}\times\frac{[2n+1]_q!}{s^{2n+1}}q^{-\frac{(2n+1)(2n+2)}{2}}\\
=&\frac{1}{s}\sum_{n=0} ^{\infty}(-1)^nq^{-(n+1)(2n+1)}\left(\frac{a}{s}\right)^{2n+1}.
\end{split}
\end{equation}
Therefore, by \eqref{76} and \eqref{77}, we obtain the following theorem.
\begin{theorem}[Transforms of $q$-cosine and $q$-sine]
\begin{gather}
{\tilde{L}}_q\left(\cos_qat\right)=\frac{1}{s}\sum_{n=0} ^{\infty}(-1)^nq^{-n(2n+1)}\left(\frac{a}{s}\right)^{2n}, \\
{\tilde{L}}_q\left(\sin_qat\right)=\frac{1}{s}\sum_{n=0} ^{\infty}(-1)^nq^{-(n+1)(2n+1)}\left(\frac{a}{s}\right)^{2n+1}.\notag
\end{gather}
\end{theorem}
We consider $q$-Laplace transforms of the second kind for the $q$-Cosine and $q$-Sine functions of the second kind. By \eqref{74} and \eqref{75}, we get
\begin{equation}\label{78}
\begin{split}
{\tilde{L}}_q\left(\Cos_q(at)\right)=&\sum_{n=0} ^{\infty}\frac{q^{\binom{2n}{2}}(-1)^na^{2n}}{[2n]_q!}\int_0 ^{\infty} t^{2n}e_q(-st)d_qt \\
=&\sum_{n=0} ^{\infty}\frac{q^{\binom{2n}{2}}(-1)^na^{2n}}{[2n]_q!}\times \frac{[2n]_q!}{s^{2n+1}}q^{-\binom{2n+1}{2}}\\
=&\frac{1}{s}\sum_{n=0} ^{\infty}\left(\frac{a}{qs}\right)^{2n}(-1)^n=\frac{q^2s}{(qs)^2+a^2},
\end{split}
\end{equation}
and
\begin{equation}\label{79}
\begin{split}
{\tilde{L}}_q\left(\Sin_q(at)\right)=&\frac{1}{2i}\left\{{\tilde{L}}_q\left(E_q(iat)\right)-{\tilde{L}}_q\left(E_q(-iat)\right)\right\}\\
=&\frac{1}{2i}\left\{\frac{q}{qs-ia}-\frac{q}{qs+ia}\right\}=\frac{1}{2i}\times\frac{2qia}{(qs)^2+a^2}=\frac{qa}{(qs)^2+a^2}.
\end{split}
\end{equation}
Therefore, by \eqref{78} and \eqref{79}, we obtain the following theorem.
\begin{theorem}[Transform of $q$-Sine and $q$-Cosine function]\label{thm11}
\begin{equation*}
{\tilde{L}}_q\left(\Cos_q(at)\right)=\frac{q^2s}{(qs)^2+a^2},~{\tilde{L}}_q\left(\Sin_q(at)\right)=\frac{qa}{(qs)^2+a^2}.
\end{equation*}
\end{theorem}
Here, we compute the $q$-Laplace transform of derivative by using the definition of the $q$-Laplace transform of the second kind. It is easy to show that
\begin{equation}\label{80}
\int_0 ^x f^{'}(t)g(t)d_qt=\left(f(x)g(t)-f(0)g(0)\right)-\int_0 ^x f(qt)g^{'}(t)d_qt,
\end{equation}
where $f^{'}(t)=D_qf(t)=\frac{\partial_qf(t)}{\partial_qt}$.

By \eqref{2}, \eqref{63} and \eqref{80}, we get
\begin{equation}\label{81}
{\tilde{L}}_q\left(f^{'}(t)\right)=\int_0 ^{\infty} f^{'}(t)e_q(-st)d_qt=-f(0)-\int_0 ^{\infty}f(qt)e_q ^{'} (-st)d_qt.
\end{equation}
From \eqref{2}, we note that
\begin{equation}\label{82}
\begin{split}
e_q ^{'}(-st)=&\frac{\partial_q}{\partial_qt}e_q(-st)=\sum_{n=0} ^{\infty}\frac{(-1)^n}{[n]_q!}s^n\frac{\partial_qt^n}{\partial_qt}=\sum_{n=1} ^{\infty}\frac{(-1)^ns^n}{[n-1]_q!}t^{n-1} \\
=&-s\sum_{n=0} ^{\infty}\frac{(-1)^n}{[n]_q!}s^nt^n=-se_q(-st).
\end{split}
\end{equation}
By \eqref{81} and \eqref{82}, we get
\begin{equation}\label{83}
\begin{split}
{\tilde{L}}_q\left(f^{'}(t)\right)=&-f(0)+s\int_0 ^{\infty}f(qt)e_q(-st)d_qt \\
=&-f(0)+sq^{-1}\int_0 ^{\infty}f(t)e_q\left(-q^{-1}st\right)d_qt\\
=&-f(0)+sq^{-1}{\tilde{F}}\left(q^{-1}s\right),
\end{split}
\end{equation}
where $f^{'}(t)=D_qf(t)=\frac{\partial_qf(t)}{\partial_qt}$. If we replace $f(t)$ by $D_qf(t)=\frac{\partial_qf(t)}{\partial_qt}$, then we have
\begin{equation}\label{84}
{\tilde{L}}_q\left(f^{(2)}(t)\right)=s^2q^{-3}{\tilde{F}}\left(q^{-2}s\right)-q^{-1}sf(0)-f^{'}(0),
\end{equation}
and
\begin{equation*}
{\tilde{L}}_q\left(f^{(3)}(t)\right)=s^3q^{-6}{\tilde{F}}\left(q^{-3}s\right)-s^2q^{-3}f(0)-sq^{-1}f^{'}(0)-f^{(2)}(0),
\end{equation*}
where $f^{(3)}(t)=D_q ^3 f(t)=\left(\frac{\partial_q}{\partial_qt}\right)^3f(t)$, $f^{(2)}(0)=\left.\left(\frac{\partial_q}{\partial_qt}\right)^2f(t)\right|_{t=0}$.
Continuing this process, we get
\begin{equation}\label{85}
{\tilde{L}}_q\left(f^{(n)}(t)\right)=s^nq^{-\binom{n+1}{2}}{\tilde{F}}\left(q^{-n}s\right)-\sum_{i=0} ^{n-1}s^{n-1-i}q^{-\binom{n-i}{2}}f^{(i)}(0),
\end{equation}
where $f^{(n)}(t)=D_q ^n f(t)=\left(\frac{\partial_q}{\partial_qt}\right)^nf(t)$, $f^{(n)}(0)=\left.\left(\frac{\partial_q}{\partial_qt}\right)^nf(t)\right|_{t=0}$. Therefore, by \eqref{85}, we obtain the following theorem.
\begin{theorem}\label{thm12}
For $n \in {\mathbb{N}}$, we have
\begin{equation*}
{\tilde{L}}_q\left(f^{(n)}(t)\right)=s^nq^{-\binom{n+1}{2}}{\tilde{F}}\left(q^{-n}s\right)-\sum_{i=0} ^{n-1}s^{n-1-i}q^{-\binom{n-i}{2}}f^{(i)}(0),
\end{equation*}
where ${\tilde{F}}(s)={\tilde{L}}_q(f(t))=\int_0 ^{\infty}e_q(-st)f(t)d_qt$, $f^{(n)}(t)=\left(\frac{\partial_q}{\partial_qt}\right)^nf(t)$.
\end{theorem}
Let us consider the following $s$-derivative for the $q$-Laplace transform of the second kind.
\begin{equation}\label{86}
\begin{split}
\frac{\partial_q{\tilde{F}}(s)}{\partial_qs}=&\frac{\partial_q{\tilde{L}}_q(f(t))}{\partial_qs}=\frac{\partial_q}{\partial_qs}\int_0 ^{\infty}e_q(-st)f(t)d_qt\\
=&\int_0 ^{\infty} \left(-te_q(-st)\right)f(t)d_qt=-{\tilde{L}}_q\left(tf(t)\right).
\end{split}
\end{equation}
Thus, by \eqref{86}, we get
\begin{equation}\label{87}
{\tilde{L}}_q\left(tf(t)\right)=-\frac{\partial_q}{\partial_qs}{\tilde{F}}(s).
\end{equation}
Continuing this process, we get
\begin{equation}\label{88}
{\tilde{L}}_q\left(t^nf(t)\right)=(-1)^n\left(\frac{\partial_q}{\partial_qs}\right)^n{\tilde{F}}(s).
\end{equation}
Therefore, by \eqref{88}, we obtain the following theorem.
\begin{theorem}\label{thm13}
For $n \in {\mathbb{N}}$, we have
\begin{equation*}
{\tilde{L}}_q\left(t^nf(t)\right)=(-1)^n\left(\frac{\partial_q}{\partial_qs}\right)^n{\tilde{F}}(s).
\end{equation*}
\end{theorem}
For example, if we take $f(t)=e_q(at)$, then
\begin{equation*}
\begin{split}
{\tilde{L}}_q\left(t^ne_q(at)\right)=&(-1)^n\left(\frac{\partial_q}{\partial_qs}\right)^n{\tilde{L}}_q\left(e_q(at)\right)=(-1)^n\left(\frac{\partial_q}{\partial_qs}\right)^n\sum_{k=0} ^{\infty}\frac{q^{-\frac{k(k+1)}{2}}}{s^{k+1}}a^k\\
=&(-1)^n\sum_{k=0} ^{\infty}a^k(-1)^kq^{-nk-\binom{k+1}{2}-\binom{n+1}{2}}\times\frac{[k+1]_q[k+2]_q\cdots [k+n]_q}{s^{k+n+1}},
\end{split}
\end{equation*}
and
\begin{equation*}
\begin{split}
{\tilde{L}}_q\left(t^nE_q(at)\right)=&(-1)^n\left(\frac{\partial_q}{\partial_qs}\right)^n{\tilde{L}}_q\left(E_q(at)\right)=(-1)^n\left(\frac{\partial_q}{\partial_qs}\right)^n\frac{q}{qs-a}\\
=&(-1)^n(-1)^nq^{n+1}[n]_q!\frac{1}{(qs-a)(q^2s-a)\cdots(q^{n+1}s-a)} \\
=&q^{n+1}\frac{[n]_q!}{(qs-a)(q^2s-a)\cdots(q^{n+1}s-a)}.
\end{split}
\end{equation*}
Finally, let us consider
\begin{equation*}
\begin{split}
L_q\left(e_q(at)f(t)\right)=&\int_0 ^{\infty}e_q(-st)e_q(at)f(t)d_qt \\
=&\sum_{n=0} ^{\infty}\frac{a^n}{[n]_q!}\int_0 ^{\infty}e_q(-st)t^nf(t)d_qt\\
=&\sum_{n=0} ^{\infty}\frac{a^n}{[n]_q!}(-1)^n{\tilde{L}}_q\left(t^nf(t)\right) \\
=&\sum_{n=0} ^{\infty}\frac{a^n}{[n]_q!}(-1)^n\left(\frac{\partial_q}{\partial_qs}\right)^n{\tilde{L}}_q(f(t))\\
=&\sum_{n=0} ^{\infty}\frac{a^n}{[n]_q!}(-1)^n\left(\frac{\partial_q}{\partial_qs}\right)^n{\tilde{F}}(s).
\end{split}
\end{equation*}


\begin{thebibliography}{10}

\bibitem {01} G. Andrews, R. Askey and R. Roy, {\it Special functions}, Cambridge University Press, Cambridge, 1999.

\bibitem {02} R. Diaz and C. Teruel, {\it $q,k$-generalized gamma and beta functions}, J. Nonlinear Math. Phys. ${\mathbf{12}}$ (2005), no. 1, 118-134.

\bibitem {03} T. Mansour and A. S. Shabani, {\it Generalization of some inequalities for the $(q_1,\cdots,q_s)$-gamma function}, Matematiche (Catania), ${\mathbf{67}}$  (2012),  o. 2, 119-130.

\bibitem {04} M. Mansour, {\it Determining the $k$-generalized gamma function $\Gamma_k(x)$ by functional equations}, Int. J. Contemp. Math. Sci., ${\mathbf{4}}$ (2009), no. 21-24, 1037-1042.

\bibitem {05}  T. Mansour, {\it Some inequalities for the $q$-gamma function}, J. Inequal. Pure Appl. Math., ${\mathbf{9}}$ (2008) no. 1, Article 18, 4 pp.

\bibitem {06}  A. de Sole and V. Kac, {\it On integral representations of $q$-gamma and $q$-beta functions},  arXiv: math. QA/0302032, 2003.

\bibitem {07}  F. Ucar and D. Albayrak, {\it On $q$-Laplace type integral operators and their applications},  J. Difference Equ. Appl., ${\mathbf{18}}$ (2012), no. 6, 1001-1014.

\bibitem {08} R. K. Yadav, S. O. Purohit and P. Nirwan, {\it On $q$-Laplace transforms of a general class of $q$-polynomials and $q$-hypergeometric functions}, Math. Maced., ${\mathbf{7}}$ (2009), 81-88.

\bibitem {09} D. G. Zill and M. R. Cullen, {\it Advanced Engineering Mathematics},  ones and Bartlett, 2005.

\end{thebibliography}
\end{document}